# Some Sufficient Conditions for the Riemann hypothesis

Choe Ryong Gil
March 6, 2012

## Abstract

The Riemann hypothesis (RH) is well known. In this paper we would show some sufficient conditions for the RH. The first condition is related with the sum of divisors function and another one is related with the Chebyshev's function.

## 1. Introduction

The function $\zeta(s)$ defined by an absolute convergent Dirichlet's series

$$\zeta(s) = \sum_{n=1}^{\infty} \frac{1}{n^s} \qquad (1)$$

in complex half-plane $\operatorname{Re} s > 1$ is called the Riemann's zeta function ([4]). The Riemann's zeta function has a simple pole with the residue 1 at $s = 1$ and except the point $s = 1$ the function $\zeta(s)$ is analytically continued to whole complex plane. And $\zeta(s)$ is expressed for $\operatorname{Re} s > 1$ as

$$\zeta(s) = \prod_{p} \left(1 - p^{-s}\right)^{-1}, \qquad (2)$$

where infinite product runs over all the prime numbers. Also for $\operatorname{Re} s > 1$ the function $\zeta(s)$ satisfies the functional equation

$$\zeta(s) = 2 \cdot (2\pi)^{s-1} \cdot \Gamma(1-s) \cdot \sin\left(\frac{\pi s}{2}\right) \cdot \zeta(1-s), \qquad (3)$$

where $\Gamma(s)$ is the gamma function ([2])

$$\Gamma(s) = \int_{0}^{+\infty} e^{-s} x^{s-1} dx. \qquad (4)$$

From the infinite product of $\zeta(s)$ the Riemann's zeta function has no zeros in $\operatorname{Re} s > 1$ and from the functional equation of $\zeta(s)$ it has trivial zeros $-2, -4, -6, \cdots$ in $\operatorname{Re} s < 0$. The zeros of $\zeta(s)$ in $0 \leq \operatorname{Re} s \leq 1$ are called the nontrivial zeros of $\zeta(s)$ ([4]). In 1859 G. Riemann conjectured that all the nontrivial zeros of $\zeta(s)$ would lie on the line $\operatorname{Re} s = 1/2$. This is just the Riemenn's hypothesis (RH). There have been published many research results on the RH. But the RH is unsolved until now ([10,11]).

To study the RH we will here consider some conditions. These conditions give us a certain possibility to prove the RH. On the bases of such new criterions, in the future, we would discuss the RH in detail.

## 2. Conditions to the sum of divisors function

In this section we will consider two sufficient conditions for the RH. Those conditions are related with the sum of divisors function.

Let $N$ be the set of natural numbers. Suppose that $p_1 = 2$, $p_2 = 3$, $p_3 = 5$, $\cdots, p_m, \cdots$ are the consecutive primes. Then $p_m$ is $m$-th prime number.

The function $\sigma(n) = \sum_{d|n} d$ is called the sum of divisors function of $n \in N$ ([2]). The relation between the function $\sigma(n)$ and the RH has an intensive expression at the Robin's inequality (RI) ([1,3,5])

$$\sigma(n) \leq e^\gamma \cdot n \cdot \log\log n, \qquad (5)$$

where $\gamma = 0.577\cdots$ is Euler's constant.

It is well known that for any $n \geq 5041$ the RI holds if and only if the RH is true. So the RI is called the Robin's criterion ([5,6,10,11]) for the RH. In the paper [1], P. Moree showed that for any odd number $n \geq 17$ the RI holds. From the infinite product of $\zeta(s)$ we put

$$R_k(p_m) = \prod_{p > p_m} \left(1 - \frac{1}{p^k}\right)^{-1}. \qquad (6)$$

The following theorem shows one property of the Riemann's zeta function for the RI. We have

**Theorem 1.** If for any $m \geq 1$ and any $k \geq 1$ it holds

$$\prod_{p \leq pm}\left(1 - \frac{1}{p_i}\right)^{-1} \cdot \prod_{p \leq pm}\left(1 - \frac{1}{p_i^{k+1}}\right) \leq e^\gamma \cdot \log p_m, \qquad (7)$$

then the RI holds.

**Proof.** Let $n = q_1^{\lambda_1} \cdot q_2^{\lambda_2} \cdots q_m^{\lambda_m}$ be a prime factorization of $n \in N$. Here $\{q_1, q_2, \cdots q_m\}$ are distinct primes and $\{\lambda_1, \lambda_2, \cdots \lambda_m\}$ are non-negative integers. We suppose $q_1 \leq q_2 \leq \cdots \leq q_m$. Let's see the theorem 1 by the induction with respect to $m$. If $m = 1$ then there exist a prime number $q$ and an integer $\lambda \geq 1$ such that $n = q^\lambda$. So if $q \geq 7$ then we have $\frac{\sigma(n)}{n} \leq 1 + \frac{1}{6} \leq 1.167$. Since $e^\gamma \cdot \log\log 7 = 1.18\cdots$, the RI holds. Suppose that if $m-1$ then the RI holds. Let's see that the RI also holds for $n = q_1^{\lambda_1} q_2^{\lambda_2} \cdots q_m^{\lambda_m}$. <u>First,</u> if $q_m \geq \log n$ then $n = q_1^{\lambda_1} q_2^{\lambda_2} \cdots q_m^{\lambda_m}$ satisfies the RI. In fact, we put $n_1 = q_1^{\lambda_1} q_2^{\lambda_2} \cdots q_{m-1}^{\lambda_{m-1}}$. Then by the assumption of the induction, the number $n_1$ satisfies the RI. On the other hand, we have

$$\frac{\sigma(n)}{n} \leq \frac{\sigma(n_1)}{n_1} \cdot \left(1 + \frac{1}{q_m - 1}\right) \leq e^\gamma \cdot \log\log n_1 \cdot \left(1 + \frac{1}{q_m - 1}\right). \qquad (8)$$

If it holds that

$$e^\gamma \cdot \log\log n_1 \cdot \left(1 + \frac{1}{q_m - 1}\right) \leq e^\gamma \cdot \log\log n, \qquad (9)$$

then $n$ satisfies the RI. In deed, it is not difficult to see that

$$\begin{aligned} q_m \cdot (\log\log n - \log\log n_1) &= \\ &= q_m \cdot \int_{\log n_1}^{\log n} \frac{dt}{t} \geq q_m \cdot \frac{(\log n - \log n_1)}{\log n} = \\ &= q_m \cdot \frac{\log q_m^{\lambda_m}}{\log n} \geq \log q_m \geq \log\log n. \end{aligned} \qquad (10)$$

Therefore $n$ satisfies the RI. <u>Next,</u> if $q_m \leq \log n$ then $n = q_1^{\lambda_1} q_2^{\lambda_2} \cdots q_m^{\lambda_m}$ satisfies the RI. In fact, put $k = \underset{1 \leq i \leq m}{Max}\{\lambda_i\}$ then we have

$$\frac{\sigma(n)}{n} = \prod_{i=1}^{m} \frac{1-q^{-\lambda_i-1}}{1-q^{-1}} \leq \prod_{i=1}^{m} \frac{1-p_i^{-\lambda_i-1}}{1-p_i^{-1}} \leq$$

$$\leq \left( \prod_{i=1}^{m} \left(1 - \frac{1}{p_i^{k+1}}\right) \right) \cdot \left( \prod_{i=1}^{m} \left(1 - \frac{1}{p_i}\right)^{-1} \right). \tag{11}$$

By the condition (7) we have

$$\log\left(\frac{\sigma(n)}{n}\right) \leq \gamma + \log\log p_m. \tag{12}$$

In general, since $p_m \leq q_m \leq \log n$, the RI holds.

This is the proof of the theorem. $\square$

**Note 1.** For real $s > 1, x > 0$, we have

$$\log R_s(x) = \sum_{p>x} \log\left(1-p^{-s}\right)^{-1} = -\sum_{p>x} \log\left(1-p^{-s}\right) =$$

$$= \sum_{p>x} \sum_{i=1}^{\infty} \frac{1}{i \cdot p^{i \cdot s}} \leq \sum_{p>x} \sum_{i=1}^{\infty} \frac{1}{p^{i \cdot s}} = \sum_{p>x} \sum_{i=1}^{\infty} \left(\frac{1}{p^s}\right)^i =$$

$$= \sum_{p>x} \frac{p^{-s}}{1-p^{-s}} = \sum_{p>x} \frac{1}{p^s - 1} \leq \frac{1}{x^{s-1}} + \sum_{n>x} \frac{1}{n^s} =$$

$$= \frac{1}{x^{s-1}} + \int_{x}^{+\infty} \frac{1}{t^s} \cdot dt = \frac{1}{x^{s-1}} + \frac{x^{1-s}}{s-1} = \frac{s}{s-1} \cdot x^{1-s}. \tag{13}$$

Hence for any $k \geq 1$ we have

$$\log\left( \prod_{i=1}^{m} \left(1 - p_i^{-k-1}\right)^{-1} \right) \to \log \zeta(k+1) \ (p_m \to \infty). \tag{14}$$

On the other hand, $\log \zeta(k+1) > 0$ and by Mertens' theorem ([4]) it holds that

$$\prod_{i=1}^{m} \left(1 - p_i^{-1}\right)^{-1} = e^{\gamma} \cdot \log p_m \left(1 + O\left(\frac{1}{\log^2 p_m}\right)\right). \tag{15}$$

And we can rewrite the condition (7) as

$$E(p_m) \le \sum_{i=1}^{m} \log\left(1 - \frac{1}{p_i^{k+1}}\right)^{-1}, \tag{16}$$

where $E(p_m) = \log \prod_{i=1}^{m}(1-p_i^{-1})^{-1} - (\log\log p_m + \gamma).$ (17)

Then from (14) and (15) we have $E(p_m) \to 0 \, (p_m \to \infty)$.

Therefore we put

$$S(k, m) = \left\{n \in N; \underset{1 \le i \le m}{Max}\{\lambda_i\} \le k\right\}, \quad S(k) = \bigcup_{m=1}^{\infty} S(k, m), \tag{18}$$

then we could know that for any $k \ge 1$ there exist only finite many numbers in the set $S(k)$ such that it doesn't satisfy the RI. This shows that The RI will hold for nearly all number except the finite numbers. □

The theorem 1 shows one sufficient condition for the RI. But below theorem shows one equivalent condition to the RI.

We have

**Theorem 2.** The RI holds if and only if it holds that

$$\limsup_{n \to \infty} \left(\frac{\sigma(n)}{n} - e^\gamma \cdot \log\log n\right) \cdot \sqrt{\log n} < +\infty. \tag{19}$$

**Proof.** Suppose that the RI holds. Then it is clear that (19) holds.

Suppose that (19) holds, but the RI doesn't hold. Then by the Robin's theorem ([5,6]), there exist constant $c > 0$, $0 < \beta < 1/2$ such that for infinitely many $n$ we have

$$e^\gamma \cdot n \cdot \log\log n + \frac{c \cdot n \cdot \log\log n}{(\log n)^\beta} \le \sigma(n). \tag{20}$$

On the other hand, since (19) holds, there exists a constant $c_0 > 0$ such that for any $n$ we have

$$\sigma(n) \leq e^{\gamma} \cdot n \cdot \log\log n + \frac{c_0 \cdot n}{\sqrt{\log n}}. \tag{21}$$

From (20) and (21), for infinitely many $n$ we have

$$e^{\gamma} \cdot n \cdot \log\log n + \frac{c \cdot n \cdot \log\log n}{(\log n)^{\beta}} \leq \sigma(n) \leq$$

$$\leq e^{\gamma} \cdot n \cdot \log\log n + \frac{c_0 \cdot n}{\sqrt{\log n}}.$$

Here since $1/2 - \beta > 0$, we have

$$0 < \frac{c}{c_0} \leq \frac{(\log\log n)^{-1}}{(\log n)^{1/2-\beta}} \to 0 \ (n \to \infty).$$

But it is a contradiction. Thus the RI holds. □

**Note 2.** In his paper [9], Ramanujan showed under the RH it holds that

$$\limsup_{n \to \infty} \left( \frac{\sigma(n)}{n} - e^{\gamma} \cdot \log\log n \right) \cdot \sqrt{\log n} \leq$$
$$\leq e^{\gamma} \cdot \left( 4 - 2\sqrt{2} + \gamma - \log 4\pi \right) = -1.39 \cdots \tag{22}$$

Therefore the theorem 2 shows that the Ramanujan's formula (22) is a condition equivalent to the RH.

We indicate that another one equivalent to the RI is that there exists a constant $c_0 \geq 1$ such that, for any $n \geq 2$,

$$\sigma(n) \leq e^{\gamma} \cdot n \cdot \log\log(c_0 \cdot n). \tag{23}$$

This is easily obtained from the theorem 2. However, in the future, we would show that the condition (23) is very important for the RH.

In this connection, we recommend the inequality

$$\sigma(n) \leq e^{\gamma} \cdot n \cdot \log\log\left( c_0 \cdot n \cdot \exp\left( \sqrt{\log n} \cdot \exp\left( \sqrt{\log\log n} \right) \right) \right), \tag{24}$$

where $c_0 \geq 1$ is a constant and $n \geq 3$. This inequality (24) is weaker than (23), but stronger than the inequality

$$\sigma(n) \leq e^{\gamma} \cdot n \cdot \log\log n + \frac{c \cdot n \cdot \exp\left(\sqrt{\log\log n}\right)}{\sqrt{\log n}}, \qquad (25)$$

where $c \geq 1$ is a constant and $n \geq 3$.

From the proof of the theorem 2, we could see that the inequality (25) is equivalent to the RH. Therefore the inequality (24) is also equivalent to the RH. This shows that the inequalities (23), (24) and (25) are equivalent to each other. However, These three inequalities have a very close relation in the proof of it. In the papers [12,13.14], we have considered specifically the inequality (23) by a new idea, which is called a sigma-index of the natural number. In particular, we gave there the proof that the inequality (24) holds unconditionally. □

## 3. A Condition to the Chebyshev's function

We will consider a more sufficient condition equivalent to the RH.

Recall that $\vartheta(x) = \sum_{p \leq x} \log p$ is called the Chebyshev's function, where $x$ is the real number and $\sum$ runs on the prime numbers not exceeding a given $x$. It is known ([3,8]) that there exists a constant $a > 0$ such that, for any $x > x_0 > 0$,

$$\vartheta(x) = x + O\left(x \cdot \exp\left(-a\sqrt{\log x}\right)\right) \qquad (26)$$

holds. And it is well known that

$$\vartheta(x) = x + O\left(\sqrt{x} \cdot \log^2 x\right) \qquad (27)$$

holds if and only if the RH is true ([7]).

We will here consider a condition related with (27).

We put

$$S_m = \sum_{n=1}^{m} \left( \frac{1-(p_{n+1}-p_n)/\log p_n}{\sqrt{p_n} \cdot \log p_n} \right). \qquad (28)$$

We have

**Theorem 3.** If $\sup_m S_m < +\infty$ the for any $p_n \geq 2$ we have

$$\vartheta(p_n) = p_n + O\left(\sqrt{p_n} \cdot \log^2 p_n\right). \qquad (29)$$

**Proof.** Let's see the theorem 2 by the induction to $n$. If $p_1 = 2$ then it is clear that $\vartheta(p_1) < p_1$. Suppose that for $p_n$ the theorem 2 holds. Then there is a constant $c_n = c_n(p_n) > 0$ such that it holds that

$$\vartheta(p_n) \leq p_n + c_n(p_n)\sqrt{p_n} \cdot \log^2 p_n. \qquad (30)$$

Let's find a condition such that theorem 2 holds for $p_{n+1}$.

Since $\vartheta(p_{n+1}) = \vartheta(p_n) + \log p_{n+1}$, it is sufficient to find a constant $c_{n+1} = c_{n+1}(p_{n+1}) > 0$ such that the inequality

$$\begin{aligned} p_n + c_n(p_n) \cdot \sqrt{p_n} \cdot \log^2 p_n + \log p_{n+1} \leq \\ \leq p_{n+1} + c_{n+1}(p_{n+1}) \cdot \sqrt{p_{n+1}} \cdot \log^2 p_{n+1} \end{aligned} \qquad (31)$$

holds. If the inequality (31) holds then we have

$$\begin{aligned} c_n(p_n) \cdot \frac{\sqrt{p_n} \cdot \log^2 p_n}{\sqrt{p_{n+1}} \cdot \log^2 p_{n+1}} + \frac{\log p_{n+1} - \log p_n}{\sqrt{p_{n+1}} \cdot \log^2 p_{n+1}} + \\ + \frac{\log p_n - (p_{n+1} - p_n)}{\sqrt{p_{n+1}} \cdot \log^2 p_{n+1}} \leq c_{n+1}(p_{n+1}). \end{aligned} \qquad (32)$$

Let's see the first term in the left hand side of (32). It is easy to see that

$$\frac{\sqrt{p_n}}{\sqrt{p_{n+1}}} = 1 - \frac{\sqrt{p_{n+1}} - \sqrt{p_n}}{\sqrt{p_{n+1}}} = 1 - \frac{p_{n+1} - p_n}{\sqrt{p_{n+1}}} \cdot \frac{1}{\left(\sqrt{p_{n+1}} + \sqrt{p_n}\right)} \qquad (33)$$

and

$$\frac{\log^2 p_n}{\log^2 p_{n+1}} = 1 - \frac{2}{\log p_{n+1}}\left(\frac{p_{n+1} - p_n}{p_n}\right) + \frac{1}{\log p_{n+1}} O\left(\frac{p_{n+1} - p_n}{p_n}\right)^2. \qquad (34)$$

So we have

$$\frac{\sqrt{p_n}\log^2 p_n}{\sqrt{p_{n+1}}\log^2 p_{n+1}} = 1 - \frac{p_{n+1}-p_n}{\sqrt{p_{n+1}}\left(\sqrt{p_{n+1}}+\sqrt{p_n}\right)} - \\ -\frac{2}{\log p_{n+1}}\left(\frac{p_{n+1}-p_n}{p_n}\right) + O\left(\frac{p_{n+1}-p_n}{p_n}\right)^2. \quad (35)$$

The second term in the left hand side of (32) is

$$\frac{\log p_{n+1}-\log p_n}{\sqrt{p_{n+1}}\log^2 p_{n+1}} = \left(\frac{p_{n+1}-p_n}{p_n}\right)\cdot\frac{1}{\sqrt{p_{n+1}}\log^2 p_{n+1}} + \\ + \frac{1}{\sqrt{p_{n+1}}\log^2 p_{n+1}}\cdot O\left(\left(\frac{p_{n+1}-p_n}{p_n}\right)^2\right). \quad (36)$$

And for any $p_n$ we have

$$\left(\frac{p_{n+1}-p_n}{p_n}\right)\cdot\frac{1}{\sqrt{p_{n+1}}\log^2 p_{n+1}} < \frac{p_{n+1}-p_n}{\sqrt{p_{n+1}}\left(\sqrt{p_{n+1}}+\sqrt{p_n}\right)}. \quad (37)$$

In general, we could suppose that $c_n(p_n) \geq \left(\sqrt{p_n}\cdot\log^2 p_n\right)^{-1}$. Then we have

$$c_n(p_n)\frac{\sqrt{p_n}\log^2 p_n}{\sqrt{p_{n+1}}\log^2 p_{n+1}} + \frac{\log p_{n+1}-\log p_n}{\sqrt{p_{n+1}}\log^2 p_{n+1}} + \\ + \frac{\log p_n - (p_{n+1}-p_n)}{\sqrt{p_{n+1}}\log^2 p_{n+1}} \leq c_n(p_n) + \frac{\log p_n - (p_{n+1}-p_n)}{\sqrt{p_n}\log^2 p_n}. \quad (38)$$

Now we take $c_{n+1}(p_{n+1})$ as

$$c_{n+1}(p_{n+1}) = c_n(p_n) + \frac{\log p_n - (p_{n+1}-p_n)}{\sqrt{p_n}\log^2 p_n} \quad (39)$$

then we have

$$\vartheta(p_{n+1}) \leq p_{n+1} + c_{n+1}(p_{n+1})\sqrt{p_{n+1}}\log^2 p_{n+1}. \quad (40)$$

From (39) we continuously have

$$c_{n+1}(p_{n+1}) = c_n(p_n) + \frac{\log p_n - (p_{n+1} - p_n)}{\sqrt{p_n} \log^2 p_n} =$$
$$= c_{n-1}(p_{n-1}) + \frac{\log p_{n-1} - (p_n - p_{n-1})}{\sqrt{p_{n-1}} \log^2 p_{n-1}} + \quad (41)$$
$$+ \frac{\log p_n - (p_{n+1} - p_n)}{\sqrt{p_n} \log^2 p_n} = \cdots = \sum_{i=1}^{n} \frac{\log p_i - (p_{i+1} - p_i)}{\sqrt{p_i} \log^2 p_i}.$$

By the assumption of the theorem, it holds that

$$c_0 = \sup_m \left( \sum_{n=1}^{m} \frac{\log p_n - (p_{n+1} - p_n)}{\sqrt{p_n} \log^2 p_n} \right) < +\infty. \quad (42)$$

Thus for any $n$ we have

$$\vartheta(p_n) \leq p_n + c_0 \sqrt{p_n} \log^2 p_n. \quad (43)$$

This is the proof of the theorem. □

**Note 3.** In the proof of the theorem 3 we assume $c_n(p_n) \geq \left(\sqrt{p_n} \cdot \log^2 p_n\right)^{-1}$.

But this is not essential. Of course, $(p_{n+1} - p_n)$ is the most important term in the series

$$\sum_{n=1}^{\infty} \frac{\log p_n - (p_{n+1} - p_n)}{\sqrt{p_n} \log^2 p_n}. \quad (44)$$

It is known that $(p_{n+1} - p_n) = O\left(\sqrt{p_n} \cdot \log^2 p_n\right)$ holds under the RH.

In fact, we could say that the essential valuation of the RH is to estimate the size of $(p_{n+1} - p_n)$. The calculation by the MATLAB to $p_m \leq 10^7$ shows

$$\sum_{p_n \leq 10^7} \frac{\log p_n - (p_{n+1} - p_n)}{\sqrt{p_n} \log^2 p_n} \leq 1.231\cdots. \quad (45)$$

Therefore we are sure that the condition $\sup_m S_m < +\infty$ would be held without any assumption. □

# References

[1] Y-J. Choie, N. Lichiardopol, P. Sole, P. Moree, "On Robin's criterion for the Riemann hypothesis", J. Theor. Nombres Bord. 19, 351-366, 2007.

[2] P. Borwin, S. Choi, B. Rooney, "The Riemann Hypothesis", Springer, 2007

[3] J. Sandor, D. S. Mitrinovic, B. Crstici, "Handbook of Number theory 1", Springer, 2006.

[4] H. L. Montgomery, R. C. Vaugnan, "Multiplicative Number Theory", Cambridge, 2006.

[5] J. C. Lagarias, "An elementary problem equivalent to the Riemann hypothesis", Amer. Math. Monthly 109, 534-543, 2002.

[6] G. Robin, "Grandes valeurs de la fonction somme des diviseurs et hypothese de Rimann", Journal of Math. Pures et appl. 63, 187-213, 1984.

[7] J. L. Nicolas, "Peties valeurs de la fonction d´Euler", Journal of Number Theory 17, 375-388, 1983.

[8] J. B. Rosser, L. Schoenfeld, " Approximate formulars for some functions of prime numbers", Illinois J. Math. 6, 64-94, 1962.

[9] J. L. Nicplas, G. Robin, Highly composite numbers by srinibasa Ramanujan, The Ramanujan Journal 1, 119-153, 1997.

[10] P. Sole, M. Planat, Robin inequality for 7-free intergers, arXiv: 1012.067v1, Dec 3,2010.

[11] G. Caveney, J. L. Nicolas, J. Sondow, Robin's theorem, primes, and a new elementary reformulation of the Riemann hypothesis, arXiv: 1110.5078v1, Oct 23, 2011.

[12] Choe R. G, One Inequality related with the Robin Inequality, February 2012, viXra.org e-Print archive, viXra1202.0030, http://vixra.org/abs/1202.0030

[13] Choe R. G, An Equivalent Inequality to the Riemann Hypothesis, February 2012, viXra.org e-Print archive, viXra1202.0048, http://vixra.org/abs/1202.0048

[14] Choe R. G, The Riemann Hypothesis and the Robin inequality, February 2012, viXra.org e-Print archive, viXra1202.0055, http://vixra.org/abs/1202.0055